\documentclass[draft]{birkjour}
\usepackage{amssymb,verbatim,enumerate,ifthen}
\usepackage[mathscr]{eucal}
\usepackage[utf8]{inputenc}
\usepackage[T1]{fontenc}
\usepackage{cite}
\def\N{\mathbb{N}}
\def\R{\mathbb{R}}
\def\Q{\mathbb{Q}}
\def\Z{\mathbb{Z}}

\def\id{\mathop{\mbox{\rm id}}\nolimits}

\newtheorem{theorem}{Theorem}
\newtheorem*{theorem*}{Theorem}
\def\Thm#1#2{\ifthenelse{\equal{#1}{*}}{\begin{theorem*}#2\end{theorem*}}
             {\begin{theorem}\label{T#1}#2\end{theorem}}}
\newtheorem{Atheorem}{Theorem}

\def\thm#1{Theorem~\ref{T#1}}

\newtheorem{proposition}[theorem]{Proposition}
\newtheorem*{proposition*}{Proposition}
\def\Prp#1#2{\ifthenelse{\equal{#1}{*}}{\begin{proposition*}#2\end{proposition*}}
             {\begin{proposition}\label{P#1}#2\end{proposition}}}

\newtheorem{corollary}[theorem]{Corollary}
\newtheorem*{corollary*}{Corollary}
\def\Cor#1#2{\ifthenelse{\equal{#1}{*}}{\begin{corollary*}#2\end{corollary*}}
             {\begin{corollary}\label{C#1}#2\end{corollary}}}
\def\cor#1{Corollary~\ref{C#1}}

\newtheorem{lemma}[theorem]{Lemma}
\newtheorem*{lemma*}{Lemma}
\def\Lem#1#2{\ifthenelse{\equal{#1}{*}}{\begin{lemma*}#2\end{lemma*}}
             {\begin{lemma}\label{L#1}#2\end{lemma}}}
\def\lem#1{Lemma~\ref{L#1}}

\newtheorem{remark}[theorem]{Remark}
\newtheorem*{remark*}{Remark}
\def\Rem#1#2{\ifthenelse{\equal{#1}{*}}{\begin{remark*}\rm #2\end{remark*}}
             {\begin{remark}\label{R#1}\rm #2\end{remark}}}

\newtheorem{example}[theorem]{Example}
\newtheorem*{example*}{Example}
\def\Exa#1#2{\ifthenelse{\equal{#1}{*}}{\begin{example*}\rm #2\end{example*}}
             {\begin{example}\label{Ex#1}\rm #2\end{example}}}

\def\eq#1{{\rm(\ref{E#1})}}
\def\Eq#1#2{\ifthenelse{\equal{#1}{*}}
  {\begin{equation*}\begin{aligned}[]#2\end{aligned}\end{equation*}}
  {\begin{equation}\begin{aligned}[]\label{E#1}#2\end{aligned}\end{equation}}}

\begin{document}

%-------------------------------------------------------------------------
% editorial commands: to be inserted by the editorial office
%
%\firstpage{1} \volume{228} \Copyrightyear{2004} \DOI{003-0001}
%\seriesextraline{This is the Concrete Title of this Book\br H.E. R and S.T.C. W, Eds.}
%\firstpage{1}

\title[Additive solvability and linear independence of the solutions]
{Additive solvability and linear independence of the solutions of a system of functional equations}

\author[E. Gselmann]{Eszter Gselmann}
\address{Institute of Mathematics, University of Debrecen, 
H-4010 Debrecen, Pf.\ 12, Hungary}
\email{gselmann@science.unideb.hu}

\author[Zs. P\'ales]{Zsolt P\'ales}
\address{Institute of Mathematics, University of Debrecen, 
H-4010 Debrecen, Pf.\ 12, Hungary}
\email{pales@science.unideb.hu}

\thanks{This work was partially supported by the European Union 
and the European Social Fund through project Supercomputer, 
the national virtual lab (grant no.: T\'{A}MOP-4.2.2.C-11/1/KONV-2012-0010).
This research was also supported by the Hungarian Scientific 
Research Fund (OTKA) Grant NK 81402.}


\date{\today}

\begin{abstract}
The aim of this paper is twofold. On one hand, the additive solvability of the system of 
functional equations
\Eq{*}{
 d_{k}(xy)=\sum_{i=0}^{k}\Gamma(i,k-i) d_{i}(x)d_{k-i}(y)
  \qquad (x,y\in \R,\,k\in\{0,\ldots,n\})
}
is studied, where $\Delta_n:=\big\{(i,j)\in\Z\times\Z\mid 0\leq i,j\mbox{ and }i+j\leq n\big\}$
and $\Gamma\colon\Delta_n\to\R$ is a symmetric function such that $\Gamma(i,j)=1$ 
whenever $i\cdot j=0$. On the other hand, the linear dependence and independence of the additive 
solutions $d_{0},d_{1},\dots,d_{n}\colon \R\to\R$ of the above system of equations is characterized. 
As a consequence of the main result, for any nonzero real derivation $d\colon\R\to\R$, the 
iterates $d^0,d^1,\dots,d^n$ of $d$ are shown to be linearly independent, 
and the graph of the mapping $x\mapsto (x,d^1(x),\dots,d^n(x))$ to be dense in $\R^{n+1}$.
\end{abstract}

%%% ----------------------------------------------------------------------
\maketitle
%%% ----------------------------------------------------------------------
\section{Introduction} 

Given a real linear space $X$, a function $a\colon\R\to X$ is called \textit{additive} if
\Eq{a}{
  a(x+y)=a(x)+a(y) \qquad(x,y\in\R).
}
It is a nontrivial fact that additive functions may satisfy further functional equations.
Among these particular additive functions the so-called derivations play an important role.
An additive function $d\colon\R\to X$ is called a \textit{derivation} (cf. \cite{Kuc85}, \cite{ZarSam58})
if it satisfies the (first-order) \textit{Leibniz Rule}:
\Eq{d}{
  d(xy)=xd(y)+yd(x) \qquad(x,y\in\R).
}
Putting $x=y=1$ into \eq{D}, we get $d(1)=0$, hence, by the $\Q$-homogeneity of additive 
functions, it follows that derivations vanish at rational numbers. Therefore, assuming that $X$ is 
equipped with a Hausdorff vector topology, the only continuous derivation is the identically zero 
function. It can be shown that derivations with weak regularity properties are necessarily continuous 
and consequently are identically equal to zero. On the other hand, there exists derivations that are
discontinuous and henceforth very irregular (see \cite{Kuc85}). More generally, for any algebraic base 
$B$ of $\R$, and for any function $d_0\colon B\to X$, there exists a unique derivation $d\colon\R\to X$ 
such that $d|_B=d_0$.

Given a real-valued derivation $d\colon \R\to\R$, one can prove by induction that the iterates 
$d^0:=\id$, $d^1:=d$, \dots, $d^{n}:=d\circ d^{n-1}$ of $d$ satisfy the following 
higher-order Leibniz Rule:
\Eq{Dk}{
 d^{k}(xy)=\sum_{i=0}^{k}\binom{k}{i}d^{i}(x)d^{k-i}(y)
  \qquad (x,y\in\R,\,k\in\{1,\ldots,n\}).
}

Motivated by this property, Heyneman--Sweedler \cite{HeySwe69} introduced the notion of $n$th-order 
derivation (in the context of functions mapping rings to modules, however, we will restrict 
ourselves only to real functions). Given $n\in \N$, a sequence of additive functions 
$d_{0},d_{1},\dots,d_{n}\colon \R\to\R$ is termed a \textit{derivation of order $n$}, if 
$d_0=\id$ and, for any $k\in\{1,\ldots,n\}$,
\Eq{dk}{
 d_{k}(xy)=\sum_{i=0}^{k}\binom{k}{i}d_{i}(x)d_{k-i}(y) \qquad (x,y\in\R)
}
is fulfilled. 

Clearly, a pair $(\id,d)$ is a first-order derivation if and only if 
$d$ is a derivation. More generally, if $d\colon\R\to\R$ is a derivation, 
then the sequence $(d^0,d^1,\dots,d^n)$ is a derivation of order $n$. However, if 
$\widetilde{d}\colon\R\to\R$ is a nontrivial derivation and $n\geq2$, then 
$(d^{0},d^{1},\dots,d^{n-1},d^{n}+\widetilde{d})$ is also an $n$th-order derivation 
where the last element is not the $n$th iterate of the derivation $d$. 

The aim of this paper is twofold. On one hand, we study the additive solvability of 
the following system of functional equations:
\Eq{dd}{
 d_{k}(xy)=\sum_{i=0}^{k}\Gamma(i,k-i) d_{i}(x)d_{k-i}(y)
  \qquad (x,y\in \R,\,k\in\{0,\ldots,n\}),
}
where
\Eq{Del}{
  \Delta_n:=\big\{(i,j)\in\Z\times\Z\mid 0\leq i,j\mbox{ and }i+j\leq n\big\},
}
and $\Gamma\colon\Delta_n\to\R$ is a symmetric function such that $\Gamma(i,j)=1$ 
whenever $i\cdot j=0$. On the other hand, we characterize the linear dependence and 
independence of the additive solutions $d_{0},d_{1},\dots,d_{n}\colon \R\to\R$ of \eq{dd}.

\section{On the additive solvability of the system of functional equations \eq{dd}}

We recall first a particular case of the following result of Ebanks \cite[Theorem 3]{Eba79b} (which 
generalizes a result of Jessen--Karpf--Thorup \cite{JesKarTho68}): 

\Lem{0}{Let $X$ be real linear space and $C,D\colon\R^2\to X$. 
Then there exists a function $f\colon\R\to X$ such that
\Eq{CD1}{
  C(x,y)&=f(x+y)-f(x)-f(y) \quad(x,y\in\R),
%\quad\mbox{and}\quad 
\\ D(x,y)&=f(xy)-xf(y)-yf(x)\quad(x,y\in\R)
}
if and only if $C,D$ satisfy the following system of equations
\Eq{CD}{
  C(x+y,z)+C(x,y)&=C(x,y+z)+C(y,z) &&(x,y,z\in\R),\\
  D(x,y)&=D(y,x) &&(x,y\in\R),\\
  D(xy,z)+zD(x,y)&=D(x,yz)+xD(y,z) &&(x,y,z\in\R), \\
  C(xz,yz)-zC(x,y)&=D(x+y,z)-D(x,z)-D(y,z) &&(x,y,z\in\R).
}}

As a trivial consequence of this result, we can characterize those 
two-variable functions that are identical to the Leibniz difference of an additive function.

\Cor{LD}{Let $X$ be a real linear space and $D\colon\R^2\to X$. Then there 
exists an additive function $f\colon\R\to X$ fulfilling functional equation
\Eq{dD}{
  D(x,y)=f(xy)-xf(y)-yf(x) \qquad(x,y\in\R)
}
if and only if $D$ satisfies
\Eq{D}{
  D(x,y)&=D(y,x) &&\qquad(x,y\in\R),\\
  D(xy,z)+zD(x,y)&=D(x,yz)+xD(y,z) &&\qquad(x,y,z\in\R),\\
  D(x+y,z)&=D(x,z)+D(y,z) &&\qquad(x,y,z\in\R).
}}

\begin{proof} Applying \lem{0} for the function $C=0$, \eq{CD1} is equivalent to the 
additivity of $f$ and \eq{dD}, and \eq{CD} reduces to \eq{D}.
\end{proof}

Our first main result offers a sufficient condition on the recursive additive solvability of the 
functional equations \eq{dd}. We deduce this result by using \cor{LD}, however, we note that another proof 
could be elaborated applyin the results of Gselmann \cite{Gse12}.

\Thm{Solv}{Let $n\geq2$ and $\Gamma\colon\Delta_n\to\R$ be a symmetric function such that 
$\Gamma(i,j)=1$ whenever $i\cdot j=0$ and
\Eq{Gam}{
  \Gamma(i+j,k)\Gamma(i,j)=\Gamma(i,j+k)\Gamma(j,k) \quad (0\leq i,j,k\text{ and }i+j+k\leq n).
}
Let $d_0=\id$ and let $d_1,\dots,d_{n-1}\colon\R\to\R$ be additive functions such that \eq{dd} holds
for $k\in\{1,\dots,n-1\}$. Then there exists an additive function $d_n\colon\R\to\R$ such that
\eq{dd} is also valid for $k=n$.}

\begin{proof} Using $\Gamma(0,n)=\Gamma(n,0)=1$, the functional equation for $d_n\colon\R\to\R$ can 
be rewritten as
\Eq{dn}{
  d_n(xy)-xd_n(y)-yd_n(x)&=D_n(x,y)
\\:=\sum_{i=1}^{n-1}&\Gamma(i,n-i) d_{i}(x)d_{k-i}(y) 
  \qquad(x,y\in\R).
}
Thus, in view of \cor{LD}, in order that there exist an additive function $d_n$ such that \eq{dn}
hold, it is necessary and sufficient that $D=D_n$ satisfy the conditions in \eq{D}. The symmetry
of $\Gamma$ implies the symmetry, the additivity of $d_1,\dots,d_{n-1}$ results the 
biadditivity of $D_n$. Thus, it suffices to prove that $D=D_n$ also satisfies the second identity
in \eq{D}. This is equivalent to showing that, for all fixed $y\in\R$, the mapping 
$(x,z)\mapsto D_n(xy,z)+zD_n(x,y)$ is symmetric. Using equations \eq{dd} for $k\in\{1,\dots,n-1\}$,
we obtain
\Eq{*}{
  D_n(xy,z)&+zD_n(x,y) \\
  &=\sum_{k=1}^{n-1}\Gamma(k,n-k) d_{k}(xy)d_{n-k}(z)
    +z\sum_{i=1}^{n-1}\Gamma(i,n-i) d_{i}(x)d_{n-i}(y) \\
  &=\sum_{k=1}^{n-1}\Gamma(k,n-k) \bigg(\sum_{i=0}^{k}\Gamma(i,k-i)d_{i}(x)d_{k-i}(y)\bigg)d_{n-k}(z)\\
    &\hspace{2cm}+z\sum_{i=1}^{n-1}\Gamma(i,n-i) d_{i}(x)d_{n-i}(y) \\
  &=\sum_{k=0}^{n}\sum_{i=0}^{k}\Gamma(k,n-k)\Gamma(i,k-i)d_{i}(x)d_{k-i}(y)d_{n-k}(z)\\
        &\hspace{2cm}-xyd_n(z)-xzd_n(y)-yzd_n(x)\\
  &=\sum_{\alpha,\beta,\gamma\geq0,\,\alpha+\beta+\gamma=n}
       \Gamma(\alpha+\beta,\gamma)\Gamma(\alpha,\beta)d_{\alpha}(x)d_{\beta}(y)d_{\gamma}(z)\\
        &\hspace{2cm}-xyd_n(z)-xzd_n(y)-yzd_n(x).
}
The sum of the last three terms in the above expression is symmetric in $(x,z)$. 
The symmetry of the first summand is the consequence of the symmetry of 
$(\alpha,\gamma)\mapsto \Gamma(\alpha+\beta,\gamma)\Gamma(\alpha,\beta)$ which 
follows from property \eq{Gam}.
\end{proof}

In what follows, we describe the nowhere zero solutions of \eq{Gam}.

\Thm{Coc}{Let $n\geq 2$ and $\Gamma\colon \Delta_{n}\to \R\setminus\{0\}$ 
be a symmetric function so that $\Gamma(i, j)=1$ whenever $i\cdot j=0$. 
Then $\Gamma$ satisfies the functional equation \eq{Gam} if and only if there exists 
a function $\gamma\colon \{0, 1, \ldots, n\}\to \R\setminus \{0\}$ such that 
\Eq{Gij}
{ \Gamma(i, j)=\frac{\gamma(i+j)}{\gamma(i)\gamma(j)} 
\qquad 
((i, j)\in \Delta_{n}).
}} 

\begin{proof}
Define the function $\gamma\colon \{0, 1, \ldots, n\}\to \R\setminus \{0\}$ through 
\Eq{*}{
 \gamma(k)= \prod_{\ell=1}^{k-1}\Gamma(\ell, 1) 
\qquad 
 (k\in\{0, 1, \ldots, n\}). 
}
The empty product being equal to $1$, we have that $\gamma(0)=\gamma(1)=1$.

To complete the proof, we have to show that, for any $(i, j)\in \Delta_n$,
\Eq{*}{
 \Gamma(i, j)=\frac{\gamma(i+j)}{\gamma(i)\gamma(j)}. 
}
This equivalent to proving that
\Eq{ij}{
 \Gamma(i,j)\prod_{\ell=1}^{i-1}\Gamma(\ell, 1) = \prod_{\ell=j}^{i+j-1}\Gamma(\ell, 1)
  \qquad((i, j)\in \Delta_n).
}
This identity trivially holds for $i=0$, $i=1$ and for any $j\in\{0,\dots,n-i\}$. 
Let $j\in \{0,\dots,n-2\}$ be fixed. We prove \eq{ij} by induction on 
$i\in\{1,\dots,n-j\}$. Assume that \eq{ij} holds for $i\in\{1,\dots,n-j-1\}$. Then,
\Eq{ij+}{
  \Gamma(i+1,j)\prod_{\ell=1}^{i}\Gamma(\ell, 1)
  &= \frac{\Gamma(i+1,j)\Gamma(i,1)}{\Gamma(i,j)}\bigg(\Gamma(i,j)\prod_{\ell=1}^{i-1}\Gamma(\ell, 1)\bigg)\\
  &= \frac{\Gamma(i+1,j)\Gamma(i,1)}{\Gamma(i,j)}\prod_{\ell=j}^{i+j-1}\Gamma(\ell, 1)\\
  &= \frac{\Gamma(i+1,j)\Gamma(i,1)}{\Gamma(i,j)\Gamma(i+j,1)}\prod_{\ell=j}^{i+j}\Gamma(\ell, 1).
}
Using \eq{Gam}, it follows that $\Gamma(i+1,j)\Gamma(i,1)=\Gamma(i,j)\Gamma(i+j,1)$, hence \eq{ij+} yields
\eq{ij} for $i+1$ instead of $i$.

Conversely, suppose that there exists a function 
$\gamma\colon \{0, 1, \ldots, n\}\to \R\setminus\{0\}$ such that 
\Eq{*}{
 \Gamma(i, j)=\frac{\gamma(i+j)}{\gamma(i)\gamma(j)} 
\qquad 
((i, j)\in \Delta_{n}). 
}
Then, for any $i, j, k\geq 0$ with $i+j+k\leq n$, we have 
\Eq{*}
{\Gamma(i+j, k)\Gamma(i, j)
&=\frac{\gamma(i+j+k)}{\gamma(i+j)\gamma(k)}\cdot \frac{\gamma(i+j)}{\gamma(i)\gamma(j)}\\
&=\frac{\gamma(i+j+k)}{\gamma(i)\gamma(j+k)}\cdot \frac{\gamma(j+k)}{\gamma(j)\gamma(k)}
 =\Gamma(i, j+k)\Gamma(j, k), 
}
which completes the proof. 
\end{proof}

When $\Gamma$ is of the form \eq{Gij}, then \thm{Solv} reduces to the following statement.

\Cor{Solv}{Let $n\geq2$ and $\gamma\colon\{0, 1, \ldots, n\}\to \R\setminus \{0\}$ with $\gamma(0)=1$.
Let $d_0=\id$ and let $d_1,\dots,d_{n-1}\colon\R\to\R$ be additive functions such that 
\Eq{DK}{
 d_{k}(xy)=\sum_{i=0}^{k}\frac{\gamma(k)}{\gamma(i)\gamma(k-i)}d_{i}(x)d_{k-i}(y) \qquad (x,y\in\R)
}
holds for $k\in\{1,\dots,n-1\}$. Then there exists an additive function $d_n\colon\R\to\R$ such that
\eq{DK} is also valid for $k=n$.}

We note that if in the above corollary $\gamma(k)=k!$, then \eq{DK} is equivalent to \eq{dk}, that is
$\id, d_1,\dots,d_{n}$ is a derivation of order $n$.

\section{A characterization of the linear dependence of additive functions}

\Thm{1}{Let $X$ be a Hausdorff locally convex linear space and let $a\colon\R\to X$ be an additive function. 
Then the following statements are equivalent:
\begin{enumerate}[(i)]
 \item There exists a nonzero continuous linear functional $\varphi\in X^*$ such that $\varphi\circ a=0$;
 \item There exists an upper semicontinuous function $\Phi:X\to\R$ such that $\Phi\not\geq0$ and $\Phi\circ 
a\geq0$;
 \item The range of $a$ is not dense in $X$, i.e., $\overline{a(\R)}\neq X$. 
\end{enumerate}}

\begin{proof}
The implication \textit{(i)$\Rightarrow$(ii)} is obvious, because $\Phi$ can be chosen as $\varphi$.

To prove \textit{(ii)$\Rightarrow$(iii)}, assume that there exists an upper semicontinuous 
function $\Phi:X\to\R$ such that $\Phi\not\geq0$ and $\Phi\circ a\geq0$. Then 
$U:=\{x\in X\mid\Phi(x)<0\}$ is a nonempty and open set. The inequality $\Phi\circ a\geq0$ 
implies that $U\cap a(\R)=\emptyset$, which proves that the range of $a$ cannot be dense in $X$. 

Finally, suppose that $\overline{a(\R)}\neq X$. By the additivity of $a$, the set $a(\R)$ is closed under
addition and multiplication by rational numbers. Therefore, the closure of $a(\R)$ is a proper closed linear 
subspace of $X$. Then, by the Hahn--Banach theorem, there exists a nonzero continuous linear functional 
$\varphi\in X^*$ which vanishes on $a(\R)$, i.e., $\varphi\circ a=0$ is satisfied.
\end{proof}

By taking $X=\R^n$, the above theorem immediately simplifies to the following consequence which 
characterizes the linear dependence of finitely many additive functions.

\Cor{1}{Let $n\in\N$ and $a_1,\dots,a_n\colon\R\to\R$ be additive functions. 
Then the following statements are equivalent:
\begin{enumerate}[(i)]
 \item The additive functions $a_1,\dots,a_n$ are linearly dependent, i.e., there exist 
 $c_1,\dots,c_n\in \R$ such that $c_1^2+\cdots+c_n^2>0$ and $c_1a_1+\dots+c_na_n=0$;
 \item There exists an upper semicontinuous function $\Phi\colon\R^n\to\R$ such that $\Phi\not\geq0$ and 
 \Eq{*}{\Phi(a_1(x),\dots,a_n(x))\geq0\qquad(x\in\R);}
 \item The set $\{(a_1(x),\dots,a_n(x))\mid x\in\R\}$ is not dense in $\R^n$. 
\end{enumerate}}

In the particular case of this corollary, namely when $\Phi$ is an indefinite quadratic form, 
the equivalence of statements \textit{(i)} and \textit{(ii)} is the main result of the paper 
\cite{Koc12} by Kocsis. A former result in this direction is due to Maksa and Rätz \cite{MakRat81P}: If two 
additive functions $a,b\colon\R\to\R$ satisfy $a(x)b(x)\geq0$ then $a$ and $b$ are linearly dependent.

\section{Linear independence of iterates of nonzero derivations}

In this section we apply \cor{1} to the particular case when the additive functions are iterates
of a real derivation. However, firstly we prove the following for higher order derivations. 

\Thm{2}{Let $n\in\N$, let $\Gamma\colon\Delta_n\to\R$ be a symmetric function such that $\Gamma(i,j)=1$
whenever $i\cdot j=0$, \eq{Gam} is satisfied and, for all $k\in\{2,\dots,n\}$ there exists 
$i\in\{1,\dots,k-1\}$
such that $\Gamma(i,k-i)\neq0$. Assume that $d_0=\id$ and $d_{1},\dots,d_{n}\colon\R\to X$ are additive 
functions satisfying \eq{dd} for all $k\in\{1,\dots,n\}$. Then the following statements are equivalent:
\begin{enumerate}[(i)]
 \item There exist $c_0,c_1,\dots,c_n\in \R$ such that $c_0^2+c_1^2+\cdots+c_n^2>0$ and
 \Eq{ddd}{c_0x+c_1d_{1}(x)+\cdots+c_nd_n(x)=0\qquad(x\in\R);}
 \item There exists an upper semicontinuous function $\Phi\colon\R^{n+1}\to\R$ such that $\Phi\not\geq0$ and 
 \Eq{*}{\Phi(x,d_{1}(x),\dots,d_{n}(x))\geq0\qquad(x\in\R);}
 \item The set $\{(x,d_{1}(x),\dots,d_{n}(x))\mid x\in\R\}$ is not dense in $\R^{n+1}$;
 \item $d_{1}=0$.
\end{enumerate}}

\begin{proof}
Applying \cor{1} to the additive functions $a_i(x)=d_{i}(x)$ $(i\in\{0,1,\dots,n\})$, it 
follows that \textit{(i)}, \textit{(ii)} and \textit{(iii)} are equivalent. The implication 
\textit{(iv)$\Rightarrow$(i)} is obvious since if $d_1=0$, then \textit{(i)} holds with
$c_1=1$ and $c_0=c_2=\cdots=c_n=0$.

Thus, it remains to show that \textit{(i)} implies \textit{(iv)}. Assume that 
\textit{(i)} holds. Then there exist a smallest $1\leq m\leq n$ and  
$c_0,\dots,c_m\in\R$ such that $c_0^2+c_1^2+\cdots+c_m^2>0$ and
\Eq{dm}{
   c_0x+c_1d_{1}(x)+\cdots+c_md_m(x)=0\qquad(x\in\R).
}
This means that the equality
\Eq{*}{
  \gamma_0x+\gamma_1d_{1}(x)+\cdots+\gamma_{m-1}d_{m-1}(x)=0\qquad(x\in\R)
}
can only hold for $\gamma_0=\cdots=\gamma_{m-1}=0$.

Observe, that $d_1(1)=\cdots=d_n(1)=0$. Indeed, $d_1(1)=0$ is a consequence of \eq{dd} when $k=1$ because
this equation means that $d_1$ is a derivation. The rest easily follows by induction on $k$ from \eq{dd}.

Putting $x=1$ into \eq{dm}, it follows that $c_0=0$. If $m=1$, then $c_1$ cannot be 
zero, hence we obtain that $d_{1}=0$.  
Thus, we may assume that the minimal $m$ for which \eq{dm} is satisfied is non-smaller than $2$. 
Replacing $x$ by $xy$ in \eq{dm} and applying \eq{dd}, for all $x,y\in\R$, we get
\Eq{*}{
   0&=\sum_{k=1}^m c_kd_k(xy)
    =\sum_{k=1}^m c_k\bigg(\sum_{i=0}^k\Gamma(i,k-i)d_{i}(x)d_{k-i}(y)\bigg) \\
   &=\sum_{k=2}^m c_k\bigg(\sum_{i=1}^{k-1}\Gamma(i,k-i)d_{i}(x)d_{k\!-\!i}(y)\bigg)
      \!+\!x\bigg(\sum_{k=1}^m c_kd_k(y)\bigg)\!+\!y\bigg(\sum_{k=1}^m c_kd_k(x)\bigg)\\
   &=\sum_{k=2}^m \sum_{i=1}^{k-1}c_k\Gamma(i,k-i)d_{i}(x)d_{k-i}(y)
    =\sum_{i=1}^{m-1} \sum_{k=i+1}^{m}c_k\Gamma(i,k-i)d_{i}(x)d_{k-i}(y) \\
   &=\sum_{i=1}^{m-1}\bigg(\sum_{j=1}^{m-i}c_{i+j}\Gamma(i,j)d_j(y)\bigg)d_{i}(x).
} 
By the minimality of $m$, it follows from the above equality that, for all $y\in\R$,
\Eq{*}{
   \sum_{j=1}^{m-i}c_{i+j}\Gamma(i,j)d_j(y)=0\qquad (i\in\{1,\dots,m-1\}).
}
Again, by the minimality of $m$, this implies that $c_{i+j}\Gamma(i,j)=0$ for $(i,j)\in\Delta_m$
with $i,j\geq1$. By the assumption of the theorem, for all $k\in\{2,\dots,n\}$ there exists 
$i\in\{1,\dots,k-1\}$ such that $\Gamma(i,k-i)\neq0$. Thus, $c_2=\dots=c_m=0$. Therefore, by \eq{dm}, 
$c_1$ cannot be equal to zero. Then \eq{dm} simplifies to $d_{1}=0$, which was to be proved.
\end{proof}

Let $n\in \N$ be arbitrary and $d\colon \R\to \R$ be a derivation. Then the $(n+1)$-tuple 
$(\id, d, d^{2},\dots,d^{n})$ is a derivation of order $n$. Thus from the 
previous theorem we immediately get the following. 

\Cor{2}{Let $n\in\N$ and let $d\colon\R\to\R$ be a derivation. 
Then the following statements are equivalent:
\begin{enumerate}[(i)]
 \item There exist $c_0,c_1,\dots,c_n\in \R$ such that $c_0^2+c_1^2+\cdots+c_n^2>0$ and
 \Eq{ddd2}{c_0x+c_1d(x)+\cdots+c_nd^n(x)=0\qquad(x\in\R);}
 \item There exists an upper semicontinuous function $\Phi\colon\R^{n+1}\to\R$ such that $\Phi\not\geq0$ and 
 \Eq{*}{\Phi(x,d(x),\dots,d^n(x))\geq0\qquad(x\in\R);}
 \item The set $\{(x,d(x),\dots,d^n(x))\mid x\in\R\}$ is not dense in $\R^{n+1}$;
 \item $d=0$.
\end{enumerate}}

\def\cprime{$'$} \def\R{\mathbb R} \def\Z{\mathbb Z} \def\Q{\mathbb Q}
  \def\C{\mathbb C}

%\nocite{}
%\bibliography{publ,funcequ,der,prob}
%\bibliographystyle{plain}

% ------------------------------------------------------------------------
\end{document}